\documentclass[reqno,11pt]{amsart} 
\usepackage{amsmath}
\usepackage{amssymb}
\usepackage{amsthm}

\newcommand\NoBlackBoxes{\global\overfullrule0pt}
\NoBlackBoxes
\theoremstyle{plain}

\begin{document}

\title{Local limit theorems for smoothed \\ 
Bernoulli and other convolutions
}

\author{Sergey G. Bobkov$^{1}$}
\thanks{1) 
School of Mathematics, University of Minnesota, Minneapolis, MN 55455 USA. Research was partially supported by NSF grant DMS-1612961.}

\author{Arnaud Marsiglietti$^{2}$}
\thanks{2) 
Department of Mathematics, University of Florida, Gainesville, FL 32611 USA. \\ Corresponding author: Arnaud Marsiglietti, Email address: a.marsiglietti@ufl.edu}


\begin{abstract}
We explore an asymptotic behavior of densities of sums of independent
random variables that are convoluted with a small continuous noise.
\end{abstract}

\maketitle
\markboth{Sergey Bobkov and Arnaud Marsiglietti}{Smoothed
Bernoulli and other convolutions}

\def\theequation{\thesection.\arabic{equation}}
\def\E{{\mathbb E}}
\def\R{{\mathbb R}}
\def\C{{\mathbb C}}
\def\P{{\mathbb P}}
\def\Z{{\mathbb Z}}
\def\H{{\rm H}}
\def\Im{{\rm Im}}
\def\Tr{{\rm Tr}}

\def\k{{\kappa}}
\def\M{{\cal M}}
\def\Var{{\rm Var}}
\def\Ent{{\rm Ent}}
\def\O{{\rm Osc}_\mu}

\def\ep{\varepsilon}
\def\phi{\varphi}
\def\vp{\varphi}
\def\F{{\cal F}}
\def\L{{\cal L}}

\def\be{\begin{equation}}
\def\en{\end{equation}}
\def\bee{\begin{eqnarray*}}
\def\ene{\end{eqnarray*}}

\vskip5mm
\section{{\bf Introduction}}
\setcounter{equation}{0}

\vskip2mm
\noindent
Let $(X_n)_{n \geq 1}$ be independent Bernoulli random variables
taking the values $\pm 1$ with probability $1/2$.
Given a random variable $X$ with density $p$, let us consider
densities $p_n$ of the normalized sums
$$
Z_n = \frac{1}{\sqrt{n}}\,(X + X_1 + \dots + X_n).
$$
By the central limit theorem, $Z_n$ are convergent weakly in distribution
to the standard normal law, which means that, as $n \rightarrow \infty$,
$$
\sup_{a < b}\, \Big|\int_a^b (p_n(x) - \varphi(x))\,dx\Big| \rightarrow 0, 
\quad {\rm where} \ \ \varphi(x) = \frac{1}{\sqrt{2\pi}}\,e^{-x^2/2}.
$$
Therefore, one may wonder whether or not this property can be sharpened as 
convergence of $p_n$ to $\varphi$ in a stronger sense. This question appears 
naturally in the area of entropic limit theorems with involved problems of 
estimation of the entropy of $X$, especially in a high-dimensional setting 
(here, we however do not discuss such applications). When $X = 0$ and the $X_k$'s 
are i.i.d., a celebrated result of Gnedenko provides necessary and sufficient 
conditions for the uniform convergence of $p_n$ when these densities exist 
([G-K], [B-RR]). Here, we will see that the presence of a non-zero noise 
$X/\sqrt{n}$ in $Z_n$ may enlarge the range of applicability of local limit 
theorems. Let us focus on the possible convergence in the $L^2$-distance
$$
\|p_n - \varphi\|_2 = 
\bigg(\int_{-\infty}^\infty |p_n(x) - \varphi(x)|^2\,dx\bigg)^{1/2}
$$
and on the uniform convergence, i.e, for the $L^\infty$-norm
$\|p_n - \varphi\|_\infty$ (which is stronger than the $L^2$-convergence).
As it turns out, the answers essentially depend on some delicate
properties of the density $p$ of $X$, as may be seen from the
following characterization in terms of the characteristic function
$$
f(t) = \E\,e^{itX} = \int_{-\infty}^\infty e^{itx}\,p(x)\,dx, \qquad t \in \R.
$$

\vskip5mm
{\bf Theorem 1.1.} {\sl If 
\be
\|p_n - \varphi\|_2 \rightarrow 0 \quad {\sl as } \ \ n \rightarrow \infty,
\en
then
\be
f(\pi k) = 0 \quad {\sl for \ all} \ \ k \in \Z, \ k \neq 0.
\en
Conversely, if $\E\,|X| < \infty$, and $f'$ is square integrable, 
then the $L^2$-convergence $(1.1)$ holds under the condition $(1.2)$.
}

\vskip5mm
Under a stronger assumption on $f$, the $L^2$-convergence of densities
may be strengthened to the uniform convergence.

\vskip5mm
{\bf Theorem 1.2.} {\sl Assume that  the condition $(1.2)$ is fulfilled.
If $\E\,|X| < \infty$, and $f'$ is integrable, then the random variables 
$Z_n$ have continuous densities $p_n$ such that
\be
\sup_x\, |p_n(x) - \varphi(x)| \rightarrow 0 \quad {\sl as } \ \ 
n \rightarrow \infty.
\en
}

The square integrability assumption in Theorem 1.1 is not 
so restrictive. By Plancherel's theorem, it may be stated in terms 
of the density of $X$ as the property
$$
\int_{-\infty}^\infty x^2 p(x)^2\,dx < \infty.
$$
This holds true as long as $p$ is bounded, and $\E X^2 < \infty$. 

As for the condition (1.2), it is of a different nature and is also fulfilled 
for a certain family of characteristic functions. This family
includes, for example, $f(t) = \frac{\sin t}{t}$ which corresponds to the uniform
distribution $U$ on the interval $[-1,1]$, and more generally
$f(t) = g(t)\,\frac{\sin t}{t}$ with an arbitrary characteristic function $g$,
which means that the distribution of $X$ contains $U$ as a component.
The condition (1.2) may also be stated explicitly in terms of the density $p$,
by virtue of the Poisson summation formula. As we will see, if $p$
has a bounded total variation, (1.2) is equivalent to the property that
$$
2\, \sum_{k \in \Z} \int_{-\infty}^\infty p(2k + x)\,p(x)\,dx \, = \, 1.
$$

As a relatively large subfamily, one may involve all characteristic functions 
$f$ that are supported on $[-\pi,\pi]$, in which case we obtain the uniform 
convergence (1.3). But, staying in a similar class, one may remove 
the assumption that $X_n$ have a Bernoulli distribution and allow
a multidimensional setting. In the sequel, we use the standard notations
$\left<\cdot,\cdot\right>$ and $|\cdot |$ to denote respectively the canonical 
inner product and the Euclidean norm in $\R^d$. A random vector 
$Y = (Y_1,\dots,Y_d)$ in $\R^d$ is said to have an isotropic distribution, if 
$$
\E \left<Y,\theta\right>^2 = |\theta|^2 \quad {\rm for \ all} \ \ \theta \in \R^d.
$$ 
Equivalently, $\E\,Y_j Y_k = \delta_{jk}$ for all $j,k \leq d$,
where $\delta_{jk}$ is the Kronecker symbol.

In the next statement, we assume that $X$ is a random vector in $\R^d$ with
characteristic function $f(t) = \E\,e^{i\left<t,X\right>}$, $t \in \R^d$, and
that $(X_n)_{n \geq 1}$ are mean zero, independent, identically distributed random 
vectors in $\R^d$ with an isotropic distribution. By the central limit theorem,
the normalized sums $Z_n$ are convergent weakly in distribution
to the standard normal law in $\R^d$ with density
\be
\varphi(x) = \frac{1}{(2\pi)^{d/2}}\,e^{-|x|^2/2}, \qquad x \in \R^d.
\en

\vskip5mm
{\bf Theorem 1.3.} {\sl
There exists $T>0$ depending on the distribution of $X_1$ with the following
property. If $f$ is supported on the ball $|t| \leq T$, then the random vectors 
$Z_n$ have continuous densities $p_n$ such that $(1.3)$ holds true. If 
$$
\beta_3 = \sup_{|\theta| = 1} \E\, |\left<X_1,\theta\right>|^3
$$ 
is finite, one may take $T = 1/\beta_3$. If $X_1$ has a non-lattice distribution,
$T$ may be arbitrary.
}

\vskip5mm
Theorems 1.1-1.2 also admit multidimensional extensions, which we discuss
in Sections 2-3. Theorem 1.3 is proved in Section 4. In Sections 5-6 
we recall the Poisson formula, including the multidimensional case, and 
discuss its applications to (1.2). In the last Section 7, we consider 
an asymptotic behavior of densities $p_n$ in dimension one without 
the property (1.2). Under mild regularity assumptions on the distribution 
of $X$, it will be shown in particular that uniformly over all $x$
$$
p_n(x) = A_n(x) \varphi(x) + O\Big(\frac{\log n}{\sqrt{n}}\Big), \qquad
A_n(x) \, = \, 2\sum_{m \in \Z} p(2m + x \sqrt{n} + n).
$$
This asymptotic representation illustrates a strong oscillatory behavior of 
the densities $p_n(x)$ for all points $x \neq 0$, which may actually
be different for even $n$ versus odd values of $n$.

\vskip5mm
\section{{\bf Multidimensional variant of Theorem 1.1}}
\setcounter{equation}{0}

\vskip2mm
\noindent
We  denote by $L^r$, $r \geq 1$, the space of all (complex-valued) functions 
$u$ on $\R^d$ with finite norm
$$
\|u\|_r = 
\bigg(\int_{\R^d} |u(x)|^r\,dx\bigg)^{1/r}.
$$

Turning to the multidimensional variant of Theorem 1.1, suppose that 
$(X_n)_{n \geq 1}$ are independent random vectors uniformly distributed 
in the discrete cube $\{-1,1\}^d$, so that their components  (coordinates) 
represent independent Bernoulli random variables. Also, let $X$ 
be a random vector in $\R^d$ with characteristic function
$$
f(t) = \E\,e^{i\left<t,X\right>}, \qquad t \in \R^d.
$$
Like the one dimensional case, if $X$ has an absolutely continuous distribution, 
the normalized sums
$$
Z_n = \frac{1}{\sqrt{n}}\,(X + X_1 + \dots + X_n)
$$
have (some) densities $p_n$. In addition, the distributions of $Z_n$ are convergent 
weakly as $n \rightarrow \infty$ to the standard normal law in $\R^d$ with density 
$\varphi$ given in (1.4). We would like to strengthen this convergence 
with respect to the $L^2$-distance $\|p_n - \varphi\|_2$.

\vskip5mm
{\bf Theorem 2.1.} {\sl If $Z_n$ have densities $p_n$ such that
\be
\|p_n - \varphi\|_2 \rightarrow 0 \quad {\sl as } \ \ n \rightarrow \infty,
\en
then
\be
f(\pi k) = 0 \quad {\sl for \ all} \ \ k \in \Z^d, \ k \neq 0.
\en
Conversely, suppose that $\E\,|X| < \infty$ and
\be
\int_{\R^d} \frac{|f(t)|\,|f'(t)|}{\|t\|^{d-1}}\,dt < \infty,
\en
where $\|t\|$ denotes the distance from the point $t \in \R^d$
to the lattice $\pi \Z^d$. Then, $Z_n$ have densities $p_n$, and
the $L^2$-convergence $(2.1)$ holds true under the condition $(2.2)$.
}

\vskip5mm
The moment assumption on $X$ guarantees that the characteristic function $f$
has a continuous derivative (gradient) $f' = \nabla f$ with its Euclidean 
norm $|f'|$, so that (2.3) makes sense. This condition implies that $f$ 
is in $L^2$ as stated in Lemma 2.2 below, hence necessarily $X$ and all $Z_n$ 
have densities. In dimension one, the condition (2.3) is fulfilled as long as
$f$ and $f'$ are in $L^2$ (by Cauchy's inequality).
If $d \geq 2$, (2.3) is a bit more complicated; it is fulfilled when
$$
\sum_{k \in \Z^d}\, \max_{t \in Q_k}\, |f(t)|\,|f'(t)| \, < \,\infty,
$$
where $Q_k = Q + \pi k$, $Q = [- \frac{\pi}{2},\frac{\pi}{2}]^d$.
This is true, for example, under the decay assumptions such as
$$
|f(t)| \, \leq \, \frac{c}{((1 + |t_1|) \dots (1 + |t_d|))^\alpha}, \qquad
|f'(t)| \, \leq \, \frac{c}{((1 + |t_1|) \dots (1 + |t_d|))^\alpha},
$$
holding for all $t = (t_1,\dots,t_d) \in \R^d$ with some constants
$\alpha > \frac{1}{2}$ and $c>0$. For instance, this is the case, when
$X$ is uniformly distributed in the cube $[-1,1]^d$.

\vskip5mm
{\bf Lemma 2.2.} {\sl If the characteristic function $f$ of the random vector 
$X$ in $\R^d$ with finite first absolute moment satisfies the condition 
$(2.3)$, and $\sum_{k \in \Z^d} |f(\pi k)|^2 < \infty$, then $X$ has 
an absolutely continuous distribution with density in $L^2$. Moreover, if
$$
\int_{\R^d} \frac{|f'(t)|}{\|t\|^{d-1}}\,dt < \infty,
$$
and $\sum_{k \in \Z^d} |f(\pi k)| < \infty$, then $X$ has 
a bounded continuous density.
}

\vskip5mm
For the proof of the lemma, as well as of Theorem 2.1 and Theorem 3.1, 
we partition $\R^d$ into the cubes $Q_k = Q + \pi k$ introduced above, so that 
$\|t\| = |t - \pi k|$ for $t \in Q_k$.

\vskip5mm
{\bf Proof.} For a given $C^1$-smooth function $w$ on $\R^d$,
consider the functions $w_k(t) = w(\pi k + t)$, $k \in \Z^d$. Since
$w_k(t) \, = \, w_k(0) + \int_0^1 \left<w_k'(\xi t),t\right> d\xi$,
we have
\be
|w_k(t)| \, \leq \, |w_k(0)| + |t| \int_0^1 |w_k'(\xi t)|\, d\xi.
\en
Change of the variable $\xi t = s$ leads to
\bee
\int_{Q_k} |w(t)|\,dt - \pi^d\, |w(\pi k)|
 & \leq &
\int_Q \int_0^1 |t|\, |w_k'(\xi t)| \,dt\,d\xi \\
 & = &
\int_Q |w_k'(s)|\, |s|
\bigg[\int_{\frac{2}{\pi}\, \|s\|_\infty}^1 \xi^{-d-1} \,d\xi\bigg]\,ds 
 \ \leq \
c_d \int_Q \frac{|w_k'(s)|}{|s|^{d-1}} \, ds
\ene
with some constant $c_d$ depending on $d$ only, where 
$\|s\|_\infty = \max_k |s_k|$ for $s = (s_1,\dots,s_d) \in \R^d$.
It follows that
\be
\int_{\R^d} |w(t)|\,dt \, \leq \, \pi^d \sum_k |w(\pi k)| + c_d
\int_{\R^d} \frac{|w'(t)|}{\|t\|^{d-1}}\,dt.
\en

For the first claim of the lemma, we apply this inequality with 
$w(t) = |f(t)|^2 = f(t) f(-t)$.
It is $C^1$-smooth and satisfies $|w'(t)| \leq 2\,|f(t)|\,|f'(t)|$.
Hence, the right-hand side of (2.5) is finite, which means that $f \in L^2$.
Hence, $X$ has density in $L^2$ as well, by the Plancherel theorem. 
Choosing $w(t) = f(t)$, we obtain that $f$ is integrable, so that
the second claim follows from the inverse Fourier formula.
\qed

\vskip5mm
Before turning to the proof of Theorem 2.1, note that
the property (2.1) is equivalent to the convergence of the $L^2$-norms
\be
\|p_n\|_2 \rightarrow \|\varphi\|_2 \quad {\rm as} \ \ n \rightarrow \infty.
\en
Indeed, formally the latter is weaker than (2.1). On the other hand,
assuming (2.6) and applying the central limit theorem with weak convergence, 
we have
$$
\|p_n - \varphi\|_2^2 \, = \, 
\|p_n\|_2^2 + \|\varphi\|_2^2 - 2\, \E\,\varphi(Z_n) \, \rightarrow \,
2\,\|\varphi\|_2^2 - 2\, \E\,\varphi(Z) \, = \, 0,
$$
where $Z$ is a standard normal random vector in $\R^d$.

Now, (2.1) requires that, for all $n$ large enough, the characteristic functions
$$
f_n(t) = \E\,e^{i\left<t,Z_n\right>} =
f\Big(\frac{t}{\sqrt{n}}\Big)\, v^n\Big(\frac{t}{\sqrt{n}}\Big)
$$
belong to $L^2$, where
$$
v(t) = \cos(t_1) \dots \cos(t_d) \quad {\rm for} \ t = (t_1,\dots,t_d) \in \R^d.
$$
Thus, introducing the characteristic function 
$g(t) = e^{-|t|^2/2}$ of $Z$ and applying the Plancherel theorem,
(2.1) may be restated as the property that
\be
\|f_n\|_2^2 \, \rightarrow \, \|g\|_2^2 \, = \, \pi^{d/2}.
\en

\vskip2mm
{\bf Proof of Theorem 2.1.}

{\bf Necessity part.} To explore the latter property, consider the integrals
$$
\|f_n\|_2^2 \, = \, \int_{\R^d} |f_n(t)|^2 \,dt \, = \,
n^{d/2} \int_{\R^d} w(t)\, v^{2n}(t) \,dt, \qquad w(t) = |f(t)|^2.
$$
Using the partition of $\R^d$ as before and
the periodicity of the cosine function, we have 
\be
\|f_n\|_2^2 \, = \,
n^{d/2} \, \sum_k \int_{Q_k} w(t)\, v^{2n}(t) \,dt \, = \,
n^{d/2} \, \sum_k I_{n,k},
\en
where
$$
I_{n,k} \, = \, \int_Q w(\pi k + t)\,v^{2n}(t) \,dt.
$$

Given $\ep > 0$, choose $t_0 > 0$ small enough such that
$w(t) \geq 1 - \ep$ in $|t| \leq t_0$. We have
$$
I_{n,0} \, \geq \, (1 - \ep) \int_{|t| \leq t_0} v^{2n}(t) \,dt \, = \,
\frac{1 - \ep}{n^{d/2}} \int_{|t| \leq t_0\sqrt{n}} 
v^{2n}\Big(\frac{t}{\sqrt{n}}\Big) \,dt,
$$
implying that
$$
\liminf_{n \rightarrow \infty} \big[n^{d/2} \, I_{n,0}\big] \, \geq \,
(1 - \ep) \int_{\R^d} e^{-|t|^2}\,dt \, = \, (1 - \ep)\,\pi^{d/2}.
$$
Since $\ep>0$ was arbitrary, we get
$$
\liminf_{n \rightarrow \infty} 
\big[n^{d/2} \, I_{n,0}\big] \geq \pi^{d/2}.
$$
A similar upper bound on $\limsup$ is obvious, and we conclude that
\be
n^{d/2} \, I_{n,0} \, \rightarrow \, \pi^{d/2}.
\en

Now, suppose that (2.2) is violated for some $k \neq 0$, that is,
$w(\pi k) > 0$. By the continuity of $w$, there exist $\ep>0$ and 
$t_0 > 0$ such that $w(\pi k + t) \geq \ep$ in $|t| \leq t_0$. Hence,
$$
I_{n,k} \geq \ep \int_{|t| \leq t_0} v^{2n}(t) \,dt =
\frac{\ep}{n^{d/2}} \int_{|t| \leq t_0\sqrt{n}}
v^{2n}\Big(\frac{t}{\sqrt{n}}\Big) \,dt,
$$
implying that
$$
\liminf_{n \rightarrow \infty} 
\big[n^{d/2} \, I_{n,k}\big] \geq \ep \int_{\R^d} e^{-|t|^2}\,dt = \ep\,\pi^{d/2}.
$$
Combining this bound with (2.9), we eventually obtain in (2.8) that
$$
\liminf_{n \rightarrow \infty}\, \|f_n\|_2^2 \, \geq \, (1+\ep)\,\pi^{d/2},
$$
which contradicts to (2.7). This proves the necessity part in Theorem 2.1.

\vskip5mm
{\bf Sufficiency part.} By Lemma 2.2, the characteristic functions $f_n$
belong to $L^2$, so that the densities $p_n$ are in $L^2$ as well. 
To prove the required relation (2.7), let us return to the representation (2.8). Recalling (2.9), our task is therefore to show that 
\be
n^{d/2} \,  \sum_{k \neq 0} I_{n,k} \, \rightarrow \, 0 \qquad
(n \rightarrow \infty).
\en

To this aim, for a fixed $k \neq 0$, using
$0 \leq \cos u \leq e^{-u^2/2}$ for $|u| \leq \frac{\pi}{2}$, we have
$$
I_{n,k} \, \leq \, J_{n,k} \, = \, \int_Q w_k(t)\,e^{-n |t|^2} \,dt,
$$
where $w_k(t) = w(\pi k + t)$, $w(t) = |f(t)|^2$ as in the proof of Lemma 2.2. 
Hence, by (2.4), and changing the variable $\xi t = s$, and then 
$\xi = \sqrt{n}\,|s|\,\frac{1}{u}$, we get
\bee
J_{n,k}
 & \leq &
\int_Q \int_0^1 |t|\, |w_k'(\xi t)|\, e^{-n |t|^2} \,dt\,d\xi \\
 & \leq &
\int_Q |w_k'(s)|\, |s| \bigg[\int_0^1 \xi^{-d-1} \, 
e^{-n |s|^2/\xi^2} \,d\xi\bigg]\,ds \\
 & \leq &
n^{-d/2} \int_Q |w_k'(s)|\, |s|^{-(d-1)} \ 
\bigg[\int_{|s|\sqrt{n}}^\infty \,u^{d-1}\, e^{-u^2} \,du\bigg]\,ds\\
 & \leq &
c_d\, n^{-d/2}
\int_Q \frac{|w_k'(s)|}{|s|^{d-1}}\, e^{-n |s|^2/2}\,ds
\ene
with some constant $c_d$ depending on the dimension, only.
Performing summation over all $k \neq 0$, we get
\be
n^{d/2} \ \sum_{k \neq 0} I_{n,k} \, \leq \, c_d
\int_{\|s\|_\infty > \frac{\pi}{2}} 
\frac{|w'(s)|}{\|s\|^{d-1}}\,e^{-n\,\|s\|^2/2} \,ds.
\en
Since $|w'(s)| \leq 2\,|f(s)|\,|f'(s)|$, and recalling the assumption (2.3),
one may apply the Lebesgue dominated convergence theorem and conclude that 
the right-hand side of (2.11) tends to zero, and thus (2.7) and (2.10) hold true.

\qed

\vskip5mm
\section{{\bf Multidimensional extension of Theorem 1.2}}
\setcounter{equation}{0}

\vskip2mm
\noindent
Keeping notations and the setting of the previous section,
the multidimensional variant of Theorem 1.2 reads as follows.

\vskip5mm
{\bf Theorem 3.1.} {\sl Let $X$ be a random vector in $\R^d$ with
$\E\,|X| < \infty$ and with characteristic function $f$ such that
\be
\int_{\R^d} \frac{|f'(t)|}{\|t\|^{d-1}}\,dt < \infty,
\en
where $\|t\|$ denotes the distance from $t \in \R^d$ to the lattice $\pi \Z^d$. 
If $f(\pi k) = 0$ for all $k \in \Z^d$, $k \neq 0$, then the normalized
sums $Z_n$ have continuous densities $p_n$ such that
\be
\sup_x\, |p_n(x) - \varphi(x)| \rightarrow 0 \quad {\sl as } \ \ 
n \rightarrow \infty.
\en
}

\vskip2mm
In dimension one, (3.1) means that $|f'|$ is integrable.
If $d \geq 2$, this condition is fulfilled when
$$
\sum_{k \in \Z^d}\, \max_{t \in Q_k}\, |f'(t)| \, < \,\infty,
$$
for example, under the decay assumptions such as
$$
|f'(t)| \, \leq \, \frac{c}{((1 + |t_1|) \dots (1 + |t_d|))^\alpha}, \qquad 
t = (t_1,\dots,t_d) \in \R^d,
$$
with some constants $\alpha > 1$ and $c>0$. Note that this is not the case, 
when $X$ is uniformly distributed in the cube $[-1,1]^d$.

This claim is very similar to Theorem 2.1, and only minor modifications
should be done in the proof of the sufficiency part.

\vskip5mm
{\bf Proof.} As before, put $v(t) = \cos(t_1) \dots \cos(t_d)$.
By Lemma 2.2, the characteristic functions 
$$
f_n(t) = f\Big(\frac{t}{\sqrt{n}}\Big)\,v^n\Big(\frac{t}{\sqrt{n}}\Big)
$$ 
are integrable. Hence, $Z_n$ have continuous densities 
given by the Fourier inversion formula
\be
p_n(x) \, = \,
\frac{1}{(2\pi)^d} \int_{\R^d} e^{-i\left<t,x\right>} f_n(t)\,dt 
 \, = \, 
n^{d/2} \sum_{k \in \Z^d} I_{n,k}(x),
\en
where
\be
I_{n,k}(x) \, = \, \frac{1}{(2\pi)^d}\,
\int_{Q_k} e^{-i\left<t,x\right>\sqrt{n}}\, f(t) v^n(t)\,dt.
\en

In particular, 
$$
n^{d/2} I_{n,0}(x) \, = \, \frac{1}{(2\pi)^d}
\int_{\sqrt{n}\,Q} e^{-i\left<t,x\right>}\, 
f\Big(\frac{t}{\sqrt{n}}\Big)\, v^n\Big(\frac{t}{\sqrt{n}}\Big) \,dt.
$$
Here, one may remove $f$ from the integrand by using the bound 
$|f(\frac{t}{\sqrt{n}}) - 1| \leq \frac{|t|}{\sqrt{n}}\, \E\,|X|$.
More precisely, this may be done at the expense of an error not exceeding
in absolute value
$$
\frac{1}{\sqrt{n}}\, \E\,|X| \int_{\sqrt{n}\,Q}
|t|\, v^n\Big(\frac{t}{\sqrt{n}}\Big) \,dt \, \leq \, 
\frac{1}{\sqrt{n}}\, \E\,|X| \int_{\sqrt{n}\,Q} |t|\,e^{-|t|^2/2} \,dt
 \, \leq \,
\frac{c}{\sqrt{n}}\, \E\,|X|
$$
up to some absolute constant $c>0$. Hence
\be
n^{d/2} I_{n,0}(x) \, = \,
\frac{1}{(2\pi)^d} \int_{\sqrt{n}\,Q}
e^{-i\left<t,x\right>}\, v^n\Big(\frac{t}{\sqrt{n}}\Big) \,dt + \theta_n(x),
\en
where $\sup_x |\theta_n(x)| \rightarrow 0$ as $n \rightarrow \infty$.
One may now turn to the approximation of $v^n$ by the Gaussian function.
With some absolute constant $c>0$, we have
$$
\Big|\cos^n\Big(\frac{u}{\sqrt{n}}\Big) - e^{-u^2/2}\Big| \, \leq \,
\frac{c}{n}\, e^{-u^2/4}, \qquad |u| \leq \frac{\pi}{2},
$$
which implies
$$
\Big|v^n\Big(\frac{t}{\sqrt{n}}\Big) - e^{-|t|^2/2}\Big| \, \leq \,
\frac{cd}{n}\, e^{-|t|^2/4}, \qquad t = (t_1,\dots,t_d) \in Q.
$$
Therefore, after another replacement, (3.5) is simplified to
$$
n^{d/2} I_{n,0}(x)  \, = \, \frac{1}{(2\pi)^d} \int_{\sqrt{n}\,Q}
e^{-i\left<t,x\right>}\, e^{-|t|^2/2} \,dt + \theta_{n,1}(x) \, = \,
\varphi(x) + \theta_{n,2}(x),
$$
where $\sup_x |\theta_{n,j}(x)| \rightarrow 0$ as $n \rightarrow \infty$.
Thus, the term in (3.3) corresponding to $k = 0$ produces the desired
normal approximation, and we are left to show that
$\sum_{k \neq 0} |I_{n,k}(x)| \rightarrow 0$ as $n \rightarrow \infty$
uniformly over all $x \in \R^d$.

For a fixed $k \neq 0$, put $w_k(t) = f(\pi k + t)$. Applying again 
$0 \leq \cos u \leq e^{-u^2/2}$ for $|u| \leq \frac{\pi}{2}$ in (3.4), we have
$$
|I_{n,k}(x)| \, \leq \, J_{n,k}, \quad 
J_{n,k} \, = \, \frac{1}{(2 \pi)^d} \int_Q |w_k(t)|\,e^{-n |t|^2/2} \,dt.
$$
Using (2.4), we therefore obtain in full analogy with the derivation
from the previous section that
$$
J_{n,k} \, \leq \, c_d \, n^{-d/2}
\int_Q \frac{|w_k'(s)|}{|s|^{d-1}}\, e^{-n |s|^2/2}\,ds
$$
with some constant $c_d$ depending on the dimension, only.
Performing summation over all $k \neq 0$, we get
\be
n^{d/2} \ \sum_{k \neq 0} |I_{n,k}(x)| \, \leq \, c_d
\int_{\|s\|_\infty > \frac{\pi}{2}} 
\frac{|f'(s)|}{\|s\|^{d-1}}\,e^{-n\,\|s\|^2/2} \,ds.
\en
Finally, by (3.1), one may apply the Lebesgue dominated convergence theorem 
and conclude that the right-hand side of (3.6) tends to zero, and thus 
(3.2) holds true.

\qed

\vskip5mm
\section{{\bf Proof of Theorem 1.3}}
\setcounter{equation}{0}

\vskip2mm
\noindent
The argument is rather standard, cf. e.g. [P1-2].
Let $v(t) = \E\,e^{i\left<t,X_1\right>}$, $t \in \R^d$, be the common 
characteristic function of $X_k$'s. 
If $f$ is supported on the ball $|t| \leq T$, the characteristic functions 
$$
f_n(t) = f\Big(\frac{t}{\sqrt{n}}\Big)\,v^n\Big(\frac{t}{\sqrt{n}}\Big)
$$ 
of the normalized sums $Z_n$ are supported on the ball
$B_n$ of radius $T\sqrt{n}$. Hence, $Z_n$ have continuous densities 
given according to the Fourier inversion formula
\be
p_n(x) \, = \, \frac{1}{(2\pi)^d} \int_{\R^d} e^{-i\left<t,x\right>}
f_n(t)\,dt \, = \, \frac{1}{(2\pi)^d} \int_{B_n} e^{-i\left<t,x\right>}
f_n(t)\,dt.
\en

In order to explore an asymptotic behavior of these integrals, first note 
that one may always choose a number $T>0$ such that, for any $0 < t_0 < T$,
\be
c(t_0) \, = \max_{t_0 \leq |t| \leq T}\, |v(t)| \, < \, 1.
\en
Moreover, by the second moment assumption, 
$$
v(t) = 1 - \frac{1}{2}\,|t|^2 + \ep(t)\,|t|^2
$$
with $\ep(t) \rightarrow 0$ as $t \rightarrow 0$. Let us choose 
$t_0 \in (0,T]$ such that $|\ep(t)| \leq \frac{1}{4}$ for all $|t| \leq t_0$. 
Then $|v(t)| \leq 1 - \frac{1}{4}\,|t|^2$ in this ball, and
$$
|f_n(t)| \leq \Big(1 - \frac{1}{4n}\,|t|^2\Big)^n \leq e^{-|t|^2/4}, \qquad
|t| \leq t_0 \sqrt{n}.
$$
Combining this estimate with (4.2), we conclude that for any sequence 
$T_n \uparrow \infty$ with $T_n \leq t_0 \sqrt{n}$,
\be
\int_{T_n \leq |t| \leq T\sqrt{n}} |f_n(t)|\,dt \, \leq \,
c^n\,(2T\sqrt{n})^d \omega_d +
\int_{|t| \geq T_n} e^{-|t|^2/4}\,dt \, \rightarrow \, 0,
\en
where $c = c(t_0)$, and $\omega_d$ denotes the volume of the $d$-dimensional 
Euclidean unit ball.

Using the principal value of the logarithm, by Taylor expansion, for 
$|t| \leq t_0$ we also have
$$
\log v(t) =
\log\Big(1 - \frac{1}{2}\,|t|^2 + \ep(t)\,|t|^2\Big) \, = \,
- \frac{1}{2}\,|t|^2 + \ep_1(t)\,|t|^2
$$
with $\ep_1(t) \rightarrow 0$ as $t \rightarrow 0$. Therefore,
$$
v^n\Big(\frac{t}{\sqrt{n}}\Big) \, = \, \exp\Big\{-\frac{1}{2}\,|t|^2 + 
\ep\Big(\frac{t}{\sqrt{n}}\Big)\,|t|^2\Big\} \, \rightarrow \, 
g(t) = e^{-|t|^2/2},
$$
where the convergence is uniform in the balls $|t| \leq T_n$
such that $T_n = o(\sqrt{n})$ as $n \rightarrow \infty$. Hence,
$$
\delta_n \, = \, \sup_{|t| \leq T_n} |f_n(t) - g(t)| \, \rightarrow \, 0.
$$
Moreover, if $T_n \uparrow \infty$ sufficiently slow,
$$
\int_{|t| \leq T_n} |f_n(t) - g(t)|\,dt \, \leq \, 
\delta_n\, (2T_n)^d \, \rightarrow \, 0
$$
as $n \rightarrow \infty$. Thus, by (4.3),
\bee
\int_{Q_n} |f_n(t) - g(t)| \,dt
 & \leq &
\int_{|t| \leq T_n} |f_n(t) - g(t)|\,dt \\
 & & + \ 
\int_{T_n \leq |t| \leq T\sqrt{n}} |f_n(t)|\,dt +
\int_{T_n \leq |t| \leq T\sqrt{n}} |g(t)|\,dt \, \rightarrow \, 0.
\ene
In view of (4.1), we obtain the desired relation  (1.3), that is,
$$
|p_n(x) - \varphi(x)| \, \leq \,\frac{1}{(2\pi)^d} 
\int_{\R^d} |f_n(t) - g(t)| \,dt\, \rightarrow \, 0.
$$

If $X_1$ has a non-lattice distribution, the property (4.2) holds true
with any $T>0$, cf. [BR-R], Section 21. 
Otherwise, let us mention how one may quantify the choice 
of $T$ satisfying (4.2). If $\xi$ is a mean zero random variable with 
$\E\,|\xi|^3 < \infty$, one has (cf. e.g. [B], Lemma 15.1)
$$
|\E\,e^{ir\xi}| \, \leq \,
\exp\Big\{-\frac{r^2}{2}\,\E \xi^2 + \frac{r^3}{3}\,\E\, |\xi|^3\Big\},
\qquad r \in \R.
$$
Applying this bound with $\xi = \left<X_1,\theta\right>$, $\theta=t/|t|$, 
$r=|t|$, $t \in \R^d$, we get
$$
|v(t)| \, \leq \, 
\exp\Big\{-\frac{|t|^2}{2} + \frac{|t|^3}{3}\,\beta_3(\theta)\Big\} \, \leq \, 
\exp\Big\{-\frac{|t|^2}{2} + \frac{|t|^3}{3}\,\beta_3\Big\},
$$
where $\beta_3(\theta) = \E\, |\left<X_1,\theta\right>|^3$. If 
$|t| \leq 1/\beta_3$, the above right-hand side does not exceed
$e^{-|t|^2/6}$. Hence, $T = 1/\beta_3$ is admissible.
\qed

\vskip5mm
{\bf Remark 4.1.} One may remove the 3rd moment assumption and
take $T = \pi$ in Theorem 1.3 (in dimension one) under the following 
hypotheses about the distribution of $X_1$
(in addition to the basic moment assumptions $\E X_1 = 0$ and $\E X_1^2 = 1$):

\vskip4mm
$a)$\, The distribution of $X_1$ is symmetric about the origin;

$b)$\, $\P\{X_1 = 0\} = 0$;

$c)$\, The distribution of $X_1$ is different than the symmetric 
Bernoulli distribution on $\{-1,1\}$.

\vskip4mm
In that case, the property (4.2) still holds true. Indeed, otherwise
take the smallest $t_0 > 0$ such that $|v(t_0)| = 1$.
This implies that $X_1$ has a lattice distribution supported on
$a + h\Z$ with $h = 2\pi/t_0$ (cf. [P2], Chapter 1, Lemma 3).
Equivalently, $X_1 = a + h \xi$ for some integer-valued 
random variable $\xi$. By the assumption $a)$, necessarily 
$a = hm/2$ for some integer $m$. Adding an integer number to $\xi$,
we may assume without loss of generality that $m = 0$ or $m = 1$.

In the first case,  $X_1 = h \xi$, so that, by $b)$, $|X_1| \geq h$ and
thus $1 = \E X_1^2 \geq h^2$. Hence $t_0 = \frac{2\pi}{h} \geq 2\pi$,
implying that (4.2) holds with any $T < 2\pi$. In the second case, 
$X_1 = h (\frac{1}{2} + \xi)$, hence $|X_1| \geq \frac{1}{2}\,h$ and
thus $1 = \E X_1^2 \geq h^2/4$. Here, by $a)$, the equality is only possible
when $\xi$ takes the values 0 and 1 with probability $1/2$, which is excluded
by $c)$. Hence $h<2$ and $t_0 = \frac{2\pi}{h} > \pi$,
implying that (4.2) holds with $T = \pi$.

\vskip5mm
\section{{\bf Poisson formula}}
\setcounter{equation}{0}

\vskip2mm
\noindent
As we mentioned before, the property (1.2), needed in Theorems 1.1-1.2 and 
their multidimensional variants, may be stated 
explicitly in terms of the density of $X$. Such a reformulation is based
on the Poisson formula which we recall in this section.

Consider the Fourier transform
$$
f(t) = \int_{\R^d} e^{i \left<t,x\right>} p(x)\,dx, \qquad t \in \R^d,
$$
for a given integrable function $p:\R^d \rightarrow \C$. The Poisson formula 
indicates that, under certain mild assumptions on $p$ (or $f$), 
we have the equality
\be
\sum_{m \in \Z^d} p(m) = \sum_{k \in \Z^d} f(2\pi k).
\en
In dimension $d=1$, it is sufficient to require that $p$ be continuous and have 
a bounded total variation on the real line. In this case, the left series in (5.1) 
is absolutely convergent, while the value of the right series is understood 
as the limit of the
corresponding symmetric sums, cf. [Z], Theorem 13.5. For higher dimensions, 
(5.1) holds true as long as $p$ belongs to the Schwarz space of functions
on $\R^d$, as mentioned in [I-K], Theorem 4.5.

Let us recall a standard argument and indicate somewhat weaker conditions
in terms of $f$, enlarging the Schwarz class, but restricting ourselves 
to the case where $p$ or $f$ are real-valued and non-negative. 

\vskip5mm
{\bf Proposition 5.1.} {\sl Let $p$ be an integrable non-negative function
on $\R^d$ whose Fourier transform $f$ is also integrable and has a continuous
derivative $f' = \nabla f$ satisfying
$$
\int_{\R^d} \frac{|f'(t)|}{\|t\|^{d-1}}\,dt < \infty,
$$
where $\|t\|$ denotes the distance from the point $t$ to the lattice $2\pi \Z^d$.
Then we have the equality $(5.1)$, in which the second series
is absolutely convergent.
}

\vskip5mm
As the next proof shows, the differentiability assumption
may slightly be relaxed, assuming that $f$ is locally Lipschitz and using 
the generalized modulus of the gradient
\be
|f'(t)| = \liminf_{s \rightarrow t}\, \frac{|f(s) - f(t)|}{|s-t|}.
\en

Note that the function $p$ in Proposition 5.1 is bounded and continuous 
(which we require below), by the integrability of $f$ and by the inverse 
Fourier formula which may be written as
$$
p(x/2\pi) = \int_{\R^d} e^{i \left<t,x\right>}\,f(-2\pi t)\,dt.
$$
This formula also shows that the role of $p$ and $f$ in (5.1) may be interchanged.
In that case, Proposition 5.1 may be restated as follows.

\vskip5mm
{\bf Proposition 5.2.} {\sl Let $p$ be an integrable, locally Lipschitz function
on $\R^d$ whose Fourier transform $f$ is integrable and non-negative. Suppose that
$$
\int_{\R^d} \frac{|p'(x)|}{\|x\|^{d-1}}\,dx < \infty,
$$
where $\|x\|$ denotes the distance from the point $x$ to the lattice $\Z^d$.
Then we have the equality $(5.1)$, in which the first series
is absolutely convergent.
}

\vskip5mm
In dimension $d=1$, the above condition on $p$ just means that $p$ has bounded
total variation on the real line, and then we arrive at the usual one-dimensional 
formulation of (5.1) under an additional assumption that $f$ is non-negative.

\vskip5mm
{\bf Proof of Proposition 5.1.} 
Let us partition $\R^d$ into the cubes $Q_k = Q + 2\pi k$, 
$Q = [-\pi,\pi]^d$, $k \in \Z^d$, and apply the bound
\be
|f(2\pi k + t) - f(2\pi k)| \, \leq \, |t|
\int_0^1 |f'(2\pi k + \xi t)|\, d\xi, \qquad t \in \R^d.
\en
It holds true as long as $f$ is locally Lipschitz, with definition (5.2) of
the modulus of the gradient of $f$. Indeed, for any $x,t \in \R^d$, the
function $u(\xi) = f(x+\xi t) - f(x)$ is locally Lipschitz on the line,
and therefore it is absolutely continuous. If $u'$ is a Radon-Nikodym derivative
of $u$, it follows from (5.2) that $|u'(\xi)| \leq |t|\,|f'(x+\xi t)|$ a.e.,
while $|u(1)| \leq \int_0^1 |u'(\xi)|\,d\xi$.

Now, arguing as in the proof of Lemma 2.2, we have
\bee
\int_Q |f(2\pi k + t) - f(2\pi k)|\, dt
 & \leq &
\int_0^1 \int_Q |f'(2\pi k + \xi t)|\,|t|\, d\xi\, dt \\
 & = &
\int_Q \bigg[|f'(2\pi k + s)|\,|s|
\int_{\frac{\|s\|_\infty}{\pi}}^1 \frac{d\xi}{\xi^{d+1}} \bigg]\, ds
 \, \leq \,
c_d \int_Q \frac{|f'(2\pi k + s)|}{|s|^{d-1}}\, ds
\ene
with some constant $c_d$ depending on $d$ only.
Hence
$$
(2\pi)^d\ |f(2\pi k)| \, \leq \, \int_{Q_k} |f(t)|\, dt + 
c_d \int_{Q_k} \frac{|f'(t)|}{\|t\|^{d-1}}\, dt.
$$
The next summation over all $k$ leads to
\be
\sum_{k \in \Z^d} |f(2\pi k)| \, \leq \, \frac{1}{(2\pi)^d}
\int_{\R^d} |f(t)|\, dt + 
\frac{c_d}{(2\pi)^d} \int_{\R^d} \frac{|f'(t)|}{\|t\|^{d-1}}\, dt \, < \, \infty,
\en
so that the second series in (5.1) is absolutely convergent. 

Next, consider the periodic function
$$
P(x) = \sum_{m \in \Z^d} p(m+x), \qquad x \in \R^d.
$$
It is a.e. finite and integrable on the unit cube $K = [0,1]^d$, since
$$
\int_K \sum_{m \in \Z^d} p(m+x)\,dx = \int_{\R^d} p(x)\,dx < \infty.
$$
Therefore, $P$ admits a multiple Fourier series expansion
$\sum_{k \in \Z^d} a_k\, e^{-2\pi i\left<k,x\right>}$ with coefficients
\bee
a_k 
 & = &
\int_K e^{2\pi i \left<k,x\right>}\,P(x)\,dx
 \, = \,
\sum_{m \in \Z^d} \int_K e^{2\pi i \left<k,x\right>}\,p(x + m)\,dx \\ 
 & = &
\sum_{m \in \Z^d} \int_{K+m} e^{2\pi i \left<k,y\right>}\,p(y)\,dy \, = \, 
\int_{\R^d} e^{2\pi i \left<k,y\right>}\,p(y)\,dy \, = \, f(2\pi k).
\ene
The Fourier series is thus absolutely convergent, and
as a consequence, $P(x) = \tilde P(x)$ a.e., where
$$
\tilde P(x) = \sum_{k \in \Z^d} f(2\pi k)\, e^{-2\pi i\left<k,x\right>}.
$$
By (5.4), $\tilde P$ represents a continuous function. 
Once $P$ is finite and continuous as well, we could
conclude that $P(x) = \tilde P(x)$ for all $x \in \R^d$. But, for $x=0$, 
the latter equality becomes the Poisson formula (5.1). 

The boundedness and continuity of $P$ (needed at zero only) 
may be explored in terms of smoothness properties of $p$. Instead, 
let us apply a smoothing argument. Using the Fourier couple on the real line,
$$
w(x) = \Big(\frac{\sin(\pi x)}{\pi x}\Big)^2, \qquad 
\hat w(t) = \Big(1 - \frac{|t|}{2\pi}\Big)^+,
$$ 
the function $w_T(x) = w(x/T)$ with a parameter $T \geq 1$ has the Fourier 
transform $\hat w_T(t) = T\,\hat w(Tt)$, $x,t \in \R$. Define
$$
w_T(x) = w_T(x_1) \dots w_T(x_d), \quad x = (x_1,\dots,x_d) \in \R^d,
$$
with its Fourier transform
$$
\hat w_T(t) = \hat w_T(t_1) \dots \hat w_T(t_d), \quad t = (t_1,\dots,t_d) \in \R^d.
$$
Put $p_T(x) = p(x)\, w_T(x)$ with the corresponding periodic function
\be
P_T(x) \, = \, \sum_{m \in \Z^d} p(m+x)\,w_T(m+x).
\en
Since $p$ is bounded, the above series is absolutely convergent. 
Indeed, using  $|w(x/T)| \leq \frac{cT^2}{1 + x^2}$, $x \in \R$,
we have, for any $x = (x_1,\dots,x_d) \in \R^d$ with $\|x\|_\infty \leq 1$,
\bee
w_T(m+x) 
 & \leq &
\frac{(cT^2)^d}{(1 + (m_1 + x_1)^2) \dots (1 + (m_d + x_d)^2)} \\
 & \leq &
\frac{(cT^2)^d}{(1 + m_1^2) \dots (1 + m_d^2)}, \qquad
m = (m_1,\dots,m_d) \in \Z^d,
\ene
where $c$ denotes an absolute constant which may be different in different places.
It follows that the sum of the series in (5.5) is uniformly bounded.
Since all terms in (5.5) are continuous in $x$, we may conclude that
$P_T$ is continuous as well.

It also follows from (5.5) and the property $w_T(0) = 1$ that
\be
\lim_{T \rightarrow \infty} P_T(0) = P(0).
\en  
It is the only place where the property that $p$ is non-negative is used.
Since $0 \leq w_T \leq 1$, we have 
$\limsup_{T \rightarrow \infty} P_T(0) \leq P(0)$. On the other hand,
since $w_T(x) \rightarrow 1$, for any fixed $N \geq 1$
$$
\liminf_{T \rightarrow \infty}\, P_T(0) \, \geq \,
\liminf_{T \rightarrow \infty} \sum_{\|m\|_\infty \leq N} p(m)\,w_T(m) \, = \,
\sum_{\|m\|_\infty \leq N} p(m).
$$
Since $N$ is arbitrary, we get
$\liminf_{T \rightarrow \infty}\, P_T(0) \geq P(0)$ and thus arrive at (5.6).

Now, the Fourier transform $f_T$ of $p_T$ represents the normalized convolution
$(2\pi)^{-d}\, f * \hat w_T$, which is integrable and satisfies
$$
\int_{\R^d} \frac{|f_T'(t)|}{\|t\|^{d-1}}\, dt < \infty.
$$
The latter follows from the equality 
$f_T' = (2\pi)^{-d}\, f' * \hat w_T = (2\pi)^{-d}\, f * \hat w_T'$ together with
the bound $\int_{\R^d} \frac{|w_T'(t)|}{\|t-s\|^{d-1}}\, dt \leq C(T)$
holding true with a constant $C(T)$ independent of $s$. Thus,
$$
\sum_{k \in \Z^d} |f_T(2\pi k)| \, < \, \infty,
$$
and we obtain the Poisson formula for the smoothed functions, that is,
\be
P_T(0) = \tilde P_T(0) \equiv \sum_{k \in \Z^d} f_T(2\pi k).
\en

In order to turn to the limit in this equality, note that 
$(2\pi)^{-d} \int_{\R^d} \hat w_T(t)\,dt = w_T(0) = w(0) = 1$, so that we may
write
$$
f_T(2\pi k) - f(2\pi k) \, = \, (2\pi)^{-d}
\int_{\R^d} (f(2\pi k + t) - f(2\pi k))\,\hat w_T(t)\,dt.
$$
Hence, by (5.3),
$$
|f_T(2\pi k) - f(2\pi k)| \, \leq \, (2\pi)^{-d}
\int_0^1 \int_{\R^d} |f'(2\pi k + \xi t)|\,|t|\,\hat w_T(t)\, d\xi\, dt.
$$
Changing the variable $\xi t = s$ and using $\hat w_T(t) = 0$ for 
$\|t\|_\infty \geq 2\pi/T$, with $|\hat w_T(t)| \leq (cT)^d$ for 
$\|t\|_\infty \leq 2\pi/T$, the last double integral may be bounded by
$$
(cT)^d
\int_{\|s\|_\infty \leq \frac{2\pi}{T}} \bigg[|f'(2\pi k + s)|\,|s|
\int_{T\, \|s\|_\infty}^1 \frac{d\xi}{\xi^{d+1}} \bigg]\, ds
 \, \leq \,
c_d \int_{\|s\|_\infty \leq \frac{2\pi}{T}} 
\frac{|f'(2\pi k + s)|}{|s|^{d-1}}\, ds
$$
with some constant $c_d$ depending on $d$ only. Hence, summing over all
$k$, we get
\be
\sum_{k \in \Z^d} |f_T(2\pi k) - f(2\pi k)| \, \leq \, c_d\, (2\pi)^{-d}
\int_{R_T} \frac{|f'(t)|}{\|t\|^{d-1}}\, dt,
\en
where $R_T = \bigcup_k \big([-\frac{2\pi}{T},\frac{2\pi}{T}]^d + 2\pi k\big)$.
This region shrinks to the lattice $2\pi \Z^d$ for growing $T$, while
the integral on the right is finite, when the integration is performed
over the whole space. Therefore, by the Lebesgue dominated convergence
theorem, both sides of (5.8) tend to zero. In particular,
$\tilde P_T(0) \rightarrow \tilde P(0)$ as $T \rightarrow \infty$.
Thus, in the limit (5.7) together with (5.6) yield the desired equality
$\tilde P(0) = P(0)$.
\qed

\vskip5mm
\section{{\bf Poisson formula for convoluted densities}}
\setcounter{equation}{0}

\vskip2mm
\noindent
Let us restate once more Propositions 5.1-5.2, assuming that $f$ is the 
characteristic function of a random vector $X$ in $\R^d$.

\vskip5mm
{\bf Proposition 6.1.} {\sl Let $\E\,|X| < \infty$, and assume that
$f$ is integrable and satisfies
\be
\int_{\R^d} \frac{|f'(t)|}{\|t\|^{d-1}}\,dt < \infty,
\en
where $\|t\|$ denotes the distance from $t$ to the lattice 
$2\pi \Z^d$. Then $X$ has a bounded continuous density $p$, and we have 
the equality $(5.1)$, in which the second series is absolutely convergent.
}

\vskip5mm
Here, the moment assumption on $X$ ensures that $f$ has a continuous
derivative $f'$.

\vskip5mm
{\bf Proposition 6.2.} {\sl Let $f$ be integrable and non-negative,
and assume that the density $p$ of $X$ is locally Lipschitz and satisfies
\be
\int_{\R^d} \frac{|p'(x)|}{\|x\|^{d-1}}\,dx < \infty,
\en
where $\|x\|$ denotes the distance from $x$ to the lattice $\Z^d$.
Then we have the equality $(5.1)$, in which the sums of both series are finite.
}

\vskip5mm
By the integrability of $f$, the random vector $X$ has a bounded continuous 
density $p$ given by the inverse Fourier formula. It implies in particular 
that $p$ has a bounded continuous derivative $p'$ as soon as
$\int_{\R^d} |t|\,|f'(t)|\,dt < \infty$. The latter condition is however not
necessary. 

Recall that in dimension $d=1$, the assumptions in Proposition 6.2 may be 
weakened. It is sufficient to require that $X$ have a continuous density 
of bounded total variation (removing any hypotheses on $f$).
This requirement may be related to the properties
of the characteristic function. For example, it is sufficient to have
(cf. e.g. [B-C-G], Proposition 5.2) that
$$
\int_{-\infty}^\infty t^2\, \big(|f(t)|^2 + |f'(t)|^2\big)\,dt < \infty.
$$

In case $d \geq 2$, the assumptions (6.1) and (6.2) are respectively
fulfilled under decay bounds
$$
|f'(t)| \, \leq \, \frac{c}{((1 + |t_1|) \dots (1 + |t_d|))^\alpha}, \qquad
|p'(x)| \, \leq \, \frac{c}{((1 + |x_1|) \dots (1 + |x_d|))^\alpha},
$$
holding for all $t = (t_1,\dots,t_d) \in \R^d$ and respectively 
$x = (x_1,\dots,x_d) \in \R^d$ with some constants $\alpha > 1$ and $c>0$.
These bounds may be strengthened to
$$
|f'(t)| \leq \frac{c}{(1 + |t|)^{\alpha d}}, \qquad 
|p'(x)| \leq \frac{c}{(1 + |x|)^{\alpha d}}.
$$
The latter is fulfilled for all functions on $\R^d$ from the Schwarz space.

Let us now turn to the density description of the condition $f(\pi k) = 0$ 
for all $k \neq 0$ appearing in Theorems 1.1-1.2 and 2.1-3.1. 
It may equivalently be stated as the property
\be
\sum_{k \in \Z^d} |f(\pi k)|^2 \, = \, 1.
\en
Note that $v(t) = |f(t/2)|^2$ is non-negative and represents the characteristic 
function of the random vector $Y = (X - X')/2$, where $X'$ is an independent 
copy of $X$. If $X$ has density $p$, the density of $Y$ is given by
$$
q(x) = 2^d \int_{\R^d} p(2x + y)\,p(y)\,dy.
$$
Hence, under the corresponding regularity assumptions, the Poisson formula
(5.1) for the couple $(q,v)$ becomes
$$
\sum_{k \in \Z^d} |f(\pi k)|^2 = \sum_{k \in \Z^d} q(k) =
2^d \sum_{k \in \Z^d} \int_{\R^d} p(2k + x)\,p(x)\,dx,
$$
which is equivalent to (6.3), if and only if
\be
\sum_{k \in \Z^d} \int_{\R^d} p(2k + x)\,p(x)\,dx \, = \, 2^{-d}.
\en

Let us precise the regularity assumptions. Since $|v'(t)| \leq 2\,|f(t)|\,|f'(t)|$,
the condition (6.1) is fulfilled as long as
\be
\int_{\R^d} \frac{|f(t)|\,|f'(t)|}{\|t\|^{d-1}}\,dt < \infty,
\en
where $\|t\|$ denotes the distance from $t$ to the lattice 
$2\pi \Z^d$. Hence, from Proposition 6.1 we obtain:

\vskip5mm
{\bf Corollary 6.3.} {\sl Let $\E\,|X| < \infty$, and assume that
$f$ is square integrable and satisfies the condition $(6.5)$. 
Then $f(\pi k) = 0$ for all $k \in \Z^d$, $k \neq 0$, if and only if
the equality $(6.4)$ holds.
}

\vskip5mm
The assumption that $f \in L^2$ implies that $X$ has a square integrable
density $p$, in which case the density $q$ is continuous. Let us also note 
that the condition (6.5) is exactly the assumption (2.3) from Theorem 2.1.
Hence, under (6.5), (6.4) is equivalent to the local limit theorem (2.1),
that is, to the property
$$
\|p_n - \varphi\|_2 \rightarrow 0 \quad {\rm as} \ \ n \rightarrow \infty.
$$

One may also develop an application of Proposition 6.2 to the density $q$
(in place of $p$). Assuming that the density $p$ has a continuous derivative,
we have that $q$ has the derivative
$$
q'(x) = 2^{d+1} \int_{\R^d} p'(2x + y)\,p(y)\,dy.
$$
To weaken the assumptions, consider the one-dimensional case.
Then, the only requirement we need to meet is that $q$ is continuous
and has a bounded total variation on the real line. The continuity is met
as long as $p \in L^2$, while $\|q\|_{{\rm TV}} \leq 2\,\|p\|_{{\rm TV}}$.
Hence, we arrive at:

\vskip5mm
{\bf Corollary 6.4.} {\sl Assume that the random variable $X$ has
a density $p$ with bounded total variation.
Then $f(\pi k) = 0$ for all $k \in \Z$, $k \neq 0$, if and only if
\be
\sum_{k \in \Z} \int_{-\infty}^\infty p(2k + x)\,p(x)\,dx \, = \, \frac{1}{2}.
\en
}

\vskip5mm
\section{{\bf Asymptotic behavior of densities without condition $(1.2)$}}
\setcounter{equation}{0}

\vskip2mm
\noindent
Let us now return to the setting of Theorem 1.2, thus restricting
ourselves to dimension $d=1$. Without the condition (1.2), the densities
$p_n(x)$ of the normalized sums $Z_n$ have an oscillating character
at all points $x \neq 0$. Here we describe a typical situation, assuming
that the density $p$ of the random variable $X$ is sufficiently regular.

\vskip5mm
{\bf Theorem 7.1.} {\sl Assume that $X$ has a continuous density $p$ 
of bounded total variation, with finite second moment. If the characteristic 
function $f = \E\,e^{itX}$ and its derivatives $f'$ and $f''$ are integrable, 
then $Z_n$ have uniformly bounded densities $p_n$ satisfying uniformly 
over all $x$
\be
p_n(x) = A_n(x) \varphi(x) + O\Big(\frac{\log n}{\sqrt{n}}\Big), \qquad
\en
where
$$
A_n(x) \, = \, 
\sum_{k \in \Z}\, e^{-i\pi k\,(x \sqrt{n} + n)}\,  f(\pi k) \, = \, 
2\sum_{m \in \Z} p(2m + x \sqrt{n} + n).
$$
}

\vskip2mm
Thus, the behavior of $p_n$ might be different for $n$ even and $n$ odd.
The point $x=0$ turns out to be special, since then the oscillatory character
disappears along even and odd values of $n$ respectively. 

\vskip5mm
{\bf Corollary 7.2.} {\sl Under the same assumptions,
\bee
\lim_{n \rightarrow \infty} p_{2n}(0) 
 & = &
\frac{2}{\sqrt{2\pi}}\, \sum_{m \in \Z} p(2m), \\ 
\lim_{n \rightarrow \infty} p_{2n+1}(0) 
 & = &
\frac{2}{\sqrt{2\pi}}\, \sum_{m \in \Z} p(2m+1).
\ene
}

\vskip2mm
{\bf Proof of Theorem 7.1.} Since $f$ is integrable, the random variables 
$Z_n$ have bounded continuous densities described by the inverse Fourier 
formula
\be
p_n(x) = \frac{1}{2\pi} \int_{-\infty}^\infty e^{-itx} f_n(t)\,dt, \qquad
x \in \R,
\en
where 
$$
f_n(t) \, = \, 
f\Big(\frac{t}{\sqrt{n}}\Big)\, \cos^n\Big(\frac{t}{\sqrt{n}}\Big)
$$
are the characteristic functions of $Z_n$. As before,
let us split the integration in (7.2) into the intervals $\frac{1}{\sqrt{n}}\,Q_k$,
$Q_k = [\pi k - \frac{\pi}{2},\pi k + \frac{\pi}{2}]$, $k \in \Z$, to get the
representation
\be
p_n(x) = \frac{1}{2\pi} \sqrt{n}\ \sum_k\, (-1)^{nk}\,
e^{-i\pi k\, x \sqrt{n}}\, I_{n,k}(x),
\en
with
\begin{eqnarray}
I_{n,k}(x)
 & = &
\int_{-\frac{\pi}{2}}^{\frac{\pi}{2}} e^{-it x \sqrt{n}}\,
f(\pi k + t)\,\cos^n(t) \,dt \nonumber \\
 & = &
f(\pi k) \int_{-\frac{\pi}{2}}^{\frac{\pi}{2}} 
e^{-it x \sqrt{n}}\, \cos^n(t) \,dt \nonumber \\
 & & + \
\int_{-\frac{\pi}{2}}^{\frac{\pi}{2}} 
e^{-it x \sqrt{n}}\, (f(\pi k + t) - f(\pi k))\,\cos^n(t) \,dt.
\end{eqnarray}

Using $\theta$, $\theta_j$ to denote quantities bounded by an absolute constant,
from the asymptotic expression
\be
\cos^n(t) = e^{-nt^2/2} + \frac{\theta}{n}\,e^{-t^2/4}, \qquad 
|t| \leq \frac{\pi}{2},
\en
we obtain that
\begin{eqnarray}
\int_{-\frac{\pi}{2}}^{\frac{\pi}{2}} e^{-it x \sqrt{n}}\, \cos^n(t) \,dt
 & = &
\int_{-\frac{\pi}{2}}^{\frac{\pi}{2}} e^{-it x \sqrt{n}}\, e^{-nt^2/2} \,dt +
\frac{\theta_1}{n} \nonumber \\
 & = &
\int_{-\infty}^{\infty} e^{-it x \sqrt{n}}\, e^{-nt^2/2} \,dt +
\frac{\theta_2}{n} \ = \
\frac{\sqrt{2\pi}}{\sqrt{n}}\,e^{-x^2/2} + \frac{\theta_2}{n}.
\end{eqnarray}
This gives an asymptotic representation for the first integral in (7.4).

The second integral has a smaller order. Put 
$\ep_n = \frac{\sqrt{2\log n}}{\sqrt{n}}$ (assuming that $\ep_n \leq \frac{\pi}{2}$).
We use $0 \leq \cos t \leq e^{-t^2/2}$, $|t| \leq \frac{\pi}{2}$, so that
$\cos^n t \leq \frac{1}{n}$ for $\ep_n < |t| < \frac{\pi}{2}$.
This implies that
\bee
\int_{\ep_n}^{\frac{\pi}{2}} |f(\pi k + t) - f(\pi k)|\,
\cos^n(t) \,dt
 & \leq &
\int_{\ep_n}^{\frac{\pi}{2}} 
\bigg[ \int_0^t |f'(\pi k + s)|\, \cos^n(t) \,ds\bigg]\,dt \\ 
 & \leq &
\int_{\ep_n}^{\frac{\pi}{2}} 
\bigg[\int_0^t |f'(\pi k + s)|\, e^{-nt^2/2} \,ds\bigg]\,dt \\
 & \leq &
\frac{2}{n} \int_0^{\frac{\pi}{2}} |f'(\pi k + s)|\,ds.
\ene
With a similar bound for the interval $-\frac{\pi}{2} < t < -\ep_n$, we get
\be
\int_{\ep_n < |t| < \frac{\pi}{2}} |f(\pi k + t) - f(\pi k)|\,
\cos^n(t) \,dt \, \leq \, \frac{2}{n}
\int_{-\frac{\pi}{2}}^{\frac{\pi}{2}} |f'(\pi k + s)|\,ds.
\en

For the interval $|t| < \ep_n$, we use the Taylor integral formula
up to the quadratic form,
$$
f(\pi k + t) - f(\pi k) = f'(\pi k) t + \int_0^t (t-s) f''(\pi k + s)\,ds.
$$
By (7.5), the linear term makes a contribution
\bee
\int_{-\ep_n}^{\ep_n} e^{-it x \sqrt{n}}\, t\,\cos^n(t) \,dt
 & = &
\int_{-\ep_n}^{\ep_n} e^{-it x \sqrt{n}}\, t\,e^{-nt^2/2} \,dt +
\frac{\theta_3}{n} \\
 & = &
\frac{1}{n}
\int_{-\ep_n\sqrt{n}}^{\ep_n\sqrt{n}} e^{-is x}\, s\, e^{-s^2/2} \,ds +
\frac{\theta_3}{n} \ = \ \frac{\theta_4}{n}.
\ene
Hence, for the integral
$$
J_{n,k}(x) = \int_{|t| < \ep_n} e^{-it x \sqrt{n}}\, (f(\pi k + t) - f(\pi k))\,
\cos^n(t) \,dt,
$$
we get
\bee
|J_{n,k}(x)|
 & \leq &
\frac{\theta_4}{n}\,|f'(\pi k)| + \bigg|\int_{-\ep_n}^{\ep_n} e^{-it x \sqrt{n}}
\int_0^t (t-s) f''(\pi k + s)\, \cos^n(t) \,dt\,ds\bigg| \\ 
 & \leq &
\frac{\theta_4}{n}\,|f'(\pi k)| + \int_{-\ep_n}^{\ep_n} \int_{-|t|}^{|t|} 
(|t|-|s|)\, |f''(\pi k + s)|\,dt\,ds \\
 & \leq &
\frac{\theta_4}{n}\,|f'(\pi k)| + 
\int_{-\ep_n}^{\ep_n} (\ep_n - |s|)^2\, |f''(\pi k + s)|\,ds \\
 & \leq &
\frac{\theta_4}{n}\,|f'(\pi k)| + \frac{2\log n}{n}
\int_{-\frac{\pi}{2}}^{\frac{\pi}{2}} |f''(\pi k + s)|\,ds.
\ene
Together with (7.6)-(7.7), we thus arrive at
\bee
\sqrt{n}\ I_{n,k}(x) 
 & = &
\Big(\sqrt{2\pi}\,e^{-x^2/2} + \frac{\theta_j}{\sqrt{n}}\,|f'(\pi k)|\Big)\,f(\pi k) \\
 & & + \
\tilde \theta_j\,\frac{\log n}{\sqrt{n}} 
\int_{-\frac{\pi}{2}}^{\frac{\pi}{2}} (|f'(\pi k + s)| + |f''(\pi k + s)|)\,ds
\ene
with bounded quantities $\theta_j$ and $\tilde \theta_j$.

To perform summation over all $k \in \Z$, first note that
$\sum_k |f(\pi k)| < \infty$, as was emphasized in Proposition 6.1. 
Similarly, $\sum_k |f'(\pi k)| < \infty$, since $f''$ is integrable.
Returning to (7.3), we thus obtain that
\bee
p_n(x)
 & = &
\Big(\varphi(x) + O(1/\sqrt{n})\Big) \sum_k\, (-1)^{nk}e^{-i\pi k\, x \sqrt{n}}\,  f(\pi k) \\
 & & + \
\theta \,\frac{\log n}{\sqrt{n}} \int_{-\infty}^{\infty} (|f'(s)| + |f''(s)|)\,ds,
\ene
that is, uniformly over all $x$
$$
p_n(x) = A_n(x) \varphi(x) + O\Big(\frac{\log n}{\sqrt{n}}\Big), \qquad
A_n(x) \, = \, \sum_{k \in \Z}\, e^{-i\pi k\,(x \sqrt{n} + n)}\,  f(\pi k).
$$

Since the factors $A_n$ are uniformly bounded, so are $p_n$.
We now apply Proposition 6.1 to the random variables
$$
\xi_n = \frac{1}{2}\,(X - x \sqrt{n} - n),
$$
whose characteristic functions and densities are given by
$$
v_n(t) = \E\,e^{it\xi_n} = e^{-it\,(x \sqrt{n} + n)/2}\,  f(t/2), \qquad
q_n(y) = 2\,p(2y + x \sqrt{n} + n).
$$
With this choice we get 
$$
A_n(x) = \sum_{k \in \Z} v_n(2\pi k) = \sum_{m \in \Z} q_n(m) = 
2 \sum_{m \in \Z} p(2m + x \sqrt{n} + n).
$$
\qed

\vskip5mm
This observation implies that we cannot hope to obtain the convergence of $p_n$ 
to $\varphi$ even in $L^1$. For example, let us consider the two-sided 
exponential distribution with density $p(x) = \frac{1}{2}\,e^{-|x|}$. 
In this case, by Corollary 7.2,
\bee
\lim_{n \rightarrow \infty} p_{2n}(0)
 & = &
\frac{1}{\sqrt{2\pi}}\, \sum_m e^{-2|m|} \ = \
\frac{1}{\sqrt{2\pi}}\,\frac{e^2+1}{e^2-1} \, > \, \varphi(0), \\
\lim_{n \rightarrow \infty} p_{2n+1}(0) 
 & = &
\frac{1}{\sqrt{2\pi}}\, \sum_m e^{-|2m+1|} \, = \,
\frac{1}{\sqrt{2\pi}}\,\frac{2e}{e^2-1} \, < \, \varphi(0).
\ene
The same expressions are obtained for the values $x = 2k/\sqrt{n}$.
The function $A_n(x)$ has period $2/\sqrt{n}$.
Let $x\sqrt{n} = 2k + h$, $k \in \Z$, $0 < h < 2$. Then along even indexes $n$, 
$$
A_n(x) \, = \, 2\sum_{m \in \Z} p(2m + h) \, = \, 
\frac{e^h + e^{2-h}}{e^2-1}.
$$
The latter expression is bounded away from zero for all $h$ small enough.
Hence, according to (7.1), we have
$\int_{-\infty}^\infty |p_n(x) - \varphi(x)|\,dx \geq c > 0$.

%\newpage

\end{document}